


\magnification=\magstep1


\hsize=14.2cm
\vsize=20cm
\hoffset=-0.4cm
\voffset=0.5cm

\input amssym.def
\input amssym.tex

\font\srm=cmr8

\font\csc=cmcsc10
\font\scsc=cmcsc10 at 8pt
\font\bcsc=cmcsc10 at 12pt
\font\title=cmr12 at 15pt

\font\teneusm=eusm10
\font\seveneusm=eusm7
\font\fiveeusm=eusm5
\newfam\eusmfam
\def\eusm{\fam\eusmfam\teneusm}
\textfont\eusmfam=\teneusm
\scriptfont\eusmfam=\seveneusm
\scriptscriptfont\eusmfam=\fiveeusm

\font\teneufm=eufm10
\font\seveneufm=eufm7
\font\fiveeufm=eufm5
\newfam\eufmfam

\textfont\eufmfam=\teneufm
\scriptfont\eufmfam=\seveneufm
\scriptscriptfont\eufmfam=\fiveeufm

\def\varGamma{{\mit \Gamma}}
\def\Re{{\rm Re}\,}
\def\Im{{\rm Im}\,}
\def\txt#1{{\textstyle{#1}}}
\def\scr#1{{\scriptstyle{#1}}}
\def\r#1{{\rm #1}}
\def\B#1{{\Bbb #1}}
\def\e#1{{\eusm #1}}

\def\sgn{{\rm sgn}}

\nopagenumbers
\def\rightheadline{\hfil{\srm Mean Values of Zeta-Functions
via Representation Theory}
\hfil\tenrm\folio}
\def\leftheadline{\tenrm\folio\hfil{\scsc Y. Motohashi}\hfil}
\def\emptyheadline{}
\headline{\ifnum\pageno=1 \emptyheadline\else
\ifodd\pageno \rightheadline \else \leftheadline\fi\fi}
\vglue 0.7cm
\centerline{\title Mean Values of Zeta-Functions}
\bigskip
\centerline{\title via Representation Theory}
\vskip 1cm
\centerline{\bcsc By Yoichi Motohashi}
\vskip 1cm
\noindent 
{\bf 0.} The aim of the present article  is to reveal a structure shared by
two basic zeta-functions in their fourth power moments. It might induce one
to ponder over the possibility to go beyond.
\smallskip
To be precise, let $\varGamma$ be a discrete subgroup of a Lie
group $G$ in the framework of $\r{GL}_2$, and let  $L_\psi$ be the
$L$-function associated with a  $\varGamma$-automorphic
function $\psi$  on $G$. Then a major subject in analytic number theory is
offered by the mean value
$$
\e{M}(L_\psi,g)=\int_{-\infty}^\infty\left|L_\psi\left(\txt{1\over2}+it\right)
\right|^2g(t)dt.\eqno(0.1)
$$  
The principal issue is to establish  an explicit spectral decomposition, and
the following two cases have so far been discussed in greater detail: 
\smallskip
\item{(1)} the fourth power
moment of the Riemann zeta-function ([9][24][25][30]), 
\item{(2)} the fourth power moment of 
the Dedekind zeta-function of the Gaussian number field ([8]),
\smallskip
\noindent
which correspond, respectively, to the specifications
$$
\eqalignno{
G&=\r{PSL}_2(\B{R}),\quad \varGamma=\r{PSL}_2(\B{Z}),&(0.2)\cr
G&=\r{PSL}_2(\B{C}),\quad \varGamma=\r{PSL}_2(\B{Z}[\sqrt{-1}]).&(0.3)
}
$$
\par
\noindent
We shall see that in these cases  the spectra come from irreducible
representations of $G$ occurring in $L^2(\varGamma\backslash G)$, and that
the resulting integral transform of $g$ has a kernel
composed of Bessel functions of representations of $G$. It should be noted
that $(2)$ can readily be extended to any imaginary quadratic number field of
class number one. Other examples that share the same structure and have been
more or less worked out are 
\smallskip
\item{(3)} the mean square of  the Dedekind zeta-function
of any quadratic number field ([29]), 
\item{(4)} the fourth power moment of the Dedekind
zeta-function of any real quadratic number field with class number one
([6]),
\item{(5)} the spectral fourth power moment of all Hecke $L$-functions,
under either $(0.2)$ or $(0.3)$ ([34]),
\smallskip
\noindent
As far as our present purpose is concerned, $(1)$ is the most fundamental, and
$(2)$ comes next endorsing our conceptual view about 
$\e{M}(L_\psi,g)$ in general. The situation with $(3)$ and $(4)$ is much similar
to $(1)$, though $(3)$ requires $\varGamma$ be replaced by a Hecke
congruence subgroup, and with $(4)$ we need to move to Hilbert
modular groups. The case $(5)$ might appear different from others, but can
in fact be regarded as an extension of either $(1)$ or $(2)$, though the situation
with $(0.3)$ is still under investigation. Note that in $(2)$--$(4)$ as well as in
$(5)$ under $(0.3)$ the twist and the average with respect to
Gr\"ossencharakteren can also be taken into account. Concerning this,
an interesting argument has been developed by P. Sarnak [37], which seems to
indicate another way to view $(2)$--$(5)$. Despite this, we shall
mostly concentrate on $(1)$ and $(2)$.
\smallskip
To the list above one may wish to add, for instance,
the mean square of individual Hecke $L$-functions associated with
cusp forms, and further more any extension to the $\r{SL}_3$
environment ([10]). It appears, however,  somewhat premature to discuss these
subjects fully in the perspective described or suggested in the present work.
We plan to return to them in the near future; nevertheless, see [18][26][31][35][36]
for instance.
\medskip
\noindent
{\csc Convention:} The weight function $g$ in $(0.1)$ is 
assumed  to be even, entire,
real on $\B{R}$, and of rapid decay in any fixed horizontal strip.  
The symbol $\zeta$ is reserved for the Riemann
zeta-function, as usual. We shall often use $F$ to denote various functions
which are specified in local contexts.  Other notations are introduced where
they are needed for the first time, and will remain effective thereafter unless
otherwise stated. Also we
stress that we are  concerned with the structural aspect, and the
asymptotic study is being more or less left aside; by the same token
we shall  often skip the discussion of convergence, naturally
under the premise that no confusion be brought in. 
\medskip
\noindent
{\csc Acknowledgement:} This article is an outcome of our recent
works which we conducted either solely or jointly
with R.W. Bruggeman. We are greatly indebted to him for his invaluable cooperation. 
The text below will thus contain various excerpts from those works, with
appropriate modifications. We thank A. Ivi\'c and
M. Jutila for their constant encouragement.
\bigskip
\noindent
{\bf 1.} First we shall describe an old observation ($(2.3)$--$(2.5)$ below) 
made in [24] (see also [30, Section 4.2]), a proper exploitation of
which was only recently actuated in our joint work [9] with Bruggeman  on
$\e{M}(\zeta^2,g)$  the fourth power moment of the Riemann zeta-function. 
\smallskip
Thus we begin with an idea of H. Weyl: 
In analytic number theory we deal mainly
with sums over rational integers
$$
\sum_{n\in I}F(n),\eqno(1.1)
$$ 
where $I$ is an interval in $\B{R}$, and $F$ an arbitrary function; and 
to estimate this there is a general principle due to Weyl. A version of it is
attributed to J.G. van der Corput and built upon the following triviality:
$$
\sum_{n\in I}F(n)={1\over{M}}\sum_{n=-\infty}^\infty\sum_{m=1}^M
F(n+m)\delta_I(n+m),
\eqno(1.2)
$$ 
where $M\ge1$ is arbitrary, and $\delta_I$ the characteristic function of
$I$. If some effective inequalities are applied to the right side, then the
result does not look trivial any more but can be a sharp tool (see
[16, Chapter 2]). In fact it could yield subconvexity bounds for
$\zeta$ which have played basic r\^{o}les in
many problems in analytic number theory. For instance the bound
$$
\zeta\left(\txt{1\over2}+it\right)\ll t^{1/6}\log t,\quad t\ge2,\eqno(1.3)
$$
is relatively easy to prove with $(1.2)$, which is significantly better
than the convexity bound following from the functional equation $(1.7)$
below. 
\smallskip 
One may regard $(1.2)$ as a kind of lifting of a one
dimensional sum to a two dimensional sum. The original problem has been
transformed into the one of finding interaction among the  non-diagonal
entries. This observation leads us to another triviality: For general double
sums we have the decomposition
$$
\sum_{m,n}F(m,n)=\Big\{\sum_{m=n}+\sum_{m<n}+\sum_{m>n}\Big\}F(m,n).
\eqno(1.4)
$$ 
This is exactly the same as what F.V. Atkinson did in his important
investigation [1] on the mean square of $\zeta$. In  [30, Section 4.1], an
essentially equivalent assertion is formulated as 
$$
\eqalignno{
\e{M}(\zeta,g)&=\int_{-\infty}^\infty\left[\Re\left\{{\Gamma'\over\Gamma}
(\txt{1\over2}+it)\right\}+2\gamma_E-\log(2\pi)\right]g(t)dt
+2\pi\Re\left\{g\left(\txt{1\over2}i\right)\right\}\cr
&+4\sum_{n=1}^\infty
d(n)\int_0^\infty
(r(r+1))^{-1/2}g_c(\log(1+1/r))\cos(2\pi{n}r)dr,&(1.5)\cr }
$$
where  $\gamma_E$ is the Euler constant, $d(n)$ the number of divisors of $n$,
and
$$
g_c(x)=\int_{-\infty}^\infty{g}(t)\cos(xt)dt.\eqno(1.6)
$$ 
We should note that $(1.5)$ implies $(1.3)$;
see Ivi\'c's lecture notes [15] for the details of
various consequences of Atkinson's result. We stress also that the last infinite
sum,  which is of a spectral nature, 
corresponds to the non-diagonal parts in $(1.4)$, and that $(1.5)$ could
be regraded as a completion of a particular application of Weyl's idea to
$\zeta$. 
\smallskip
This completion was realized in [30] via the functional equation
$$
\zeta(1-s)=2(2\pi)^{-s}\cos\left(\txt{1\over2}\pi
s\right)\Gamma(s)\zeta(s).
\eqno(1.7)
$$ 
On the other hand, in [1] had been applied the Poisson summation formula
$$
\sum_{n=1}^\infty F(n)=\int_0^\infty F(r)dr+
2\sum_{n=1}^\infty\int_0^\infty F(r)\cos(2\pi nr)dr, \eqno(1.8)
$$ 
where $F$ is assumed to be smooth and of fast decay on
$(0,\infty)$. However, this difference is superficial, for $(1.7)$ and
$(1.8)$ are equivalent to each other. In fact, the
kernel function $\cos(2\pi r)$ in $(1.8)$, as well as in $(1.5)$, is related
to $(1.7)$ via the Mellin transform: 
$$
\int_0^\infty\cos(2\pi r)r^{s-1}dr
=(2\pi)^{-s}\cos\left(\txt{1\over2}\pi s\right)\Gamma(s),
\quad 0<\Re s<1,\eqno(1.9)
$$
or, more precisely, via the convolution relation
$$
\int_0^\infty F(r)\cos(2\pi nr)
dr={1\over2\pi i}\int_{(\alpha)}F^*(1-s)
\cos\left(\txt{1\over2}\pi s\right)\Gamma(s)(2\pi n)^{-s}ds,\eqno(1.10)
$$
where $F^*$ is the Mellin transform of $F$, and $(\alpha)$
the vertical line $\Re s=\alpha$ with $\alpha>0$. Summing $(1.10)$ over $n$,
and applying $(1.7)$, we are led to $(1.8)$. Reversing the reasoning, one
may reach $(1.7)$ with a combination of $(1.8)$ and $(1.9)$. 
\smallskip
In terms of the representation theory of the Lie group $\B{R}$, 
the identity $(1.8)$ is understood as a  transformation formula of $F$ that
gives a way to compute the values of projections of the Poincar\'e series
$$
\sum_{n\in\B{Z}}F(n+x),\quad x\in\B{R},\eqno(1.11)
$$ 
to irreducible subspaces of $L^2(\B{Z}\backslash\B{R})$. The appearance 
in $(1.8)$ of the Fourier-cosine transformation indicates the nature of the
harmonic analysis over the group $\B{R}$ in which $\B{Z}$ is a discrete
subgroup.
\smallskip
A motivation of the present article is to see how far extends the above
harmonic structure lying behind the important expansion $(1.5)$. 
\bigskip
\noindent
{\bf 2.} In this context, one may ponder whether 
the dissection mode in $(1.4)$ is optimal or not. It is certainly not in
general. An option to tune it up is to refine the notion of being diagonal.
A simple arithmetic way to do this is to replace $(1.4)$ by
$$
\sum_{m,n}F(m,n)=\left\{\sum_{km=ln}+\sum_{km<ln}+\sum_{km>ln}\right\}
F(m,n),\eqno(2.1)
$$ 
where $k,\,l$ are arbitrary non-zero integers. To extract information
from all of these dissections we multiply both sides by a weight
$W(k,l)$ and sum over all $k,\,l$. We get
$$
\Bigg(\sum_{k,l}W(k,l)\Bigg)\left(\sum_{m,n}F(m,n)\right)
=\left\{\sum_{km=ln}+\sum_{km<ln}+\sum_{km>ln}\right\}W(k,l)
F(m,n).\eqno(2.2)
$$ 
As before, we may regard $(2.2)$ to be a lifting of a double sum to a four
dimensional sum; and we are led to the following trivial decomposition:
$$
\sum_{k,l,m,n}F(k,l,m,n)=\left\{\sum_{km=ln}+\sum_{km<ln}+\sum_{km>ln}
\right\}F(k,l,m,n).\eqno(2.3)
$$ 
\par
Then we take a new viewpoint: We regard quadruple sums as sums over
$2\times2$ integral matrices $M$. The last identity thus becomes
$$
\sum_{M}F(M)=\left\{\sum_{|M|=0}+\sum_{|M|>0} +\sum_{|M|<0}\right\}
F(M).\eqno(2.4)
$$ 
Invoking Hecke's representatives of integral matrices
with a given determinant, we have further
$$
\sum_{|M|>0}F(M)=\sum_{f=1}^\infty\sqrt{f}\,T_fP_{F_f}(1),\quad
P_F(\r{g})=\sum_{\gamma\in\varGamma}F\left(\gamma\r{g}\right),
\quad\r{g}\in
G,\eqno(2.5)
$$ 
with $(0.2)$, where $F_f(\r{g})=F(\sqrt{f}\cdot\r{g})$ and
$$
T_f F(\r{g})={1\over\sqrt{f}}\sum_{d|f}\sum_{b=1}^d
F\left(\left[\matrix{\sqrt{f}/d&b\cr&d/\sqrt{f}
}\right]\r{g}\right).\eqno(2.6)
$$ 
\smallskip
Hereby we find a relation lying between the dissection
$(1.4)$ and the theory of $\varGamma$-auto\-morphic functions
on $G$. If $F$ is sufficiently smooth one may apply the harmonic analysis
over $\varGamma\backslash G$ to  the Poincar\'e series $P_F$. 
Thus, we are now to find the transformation formula of $F$ that gives a way to
compute the values of projections of $P_F$ to irreducible subspaces of 
$L^2(\varGamma\backslash G)$. In other words, our interest is in fixing the
kernel function corresponding  to $\cos(2\pi r)$ of
$(1.8)$, via the theory of irreducible representations of $G$ occurring in
$L^2(\varGamma\backslash G)$. In the application to 
$\e{M}(\zeta^2,g)$, as developed in [9], those four integral
variables in $(2.3)$ correspond to the four zeta-values in an obvious way. 
In the subsequent sections we shall show the salient points of the 
pertinent parts  of [9] and briefly describe a solution to the last problem.
\smallskip
It should be stressed here that the restriction
of the sum $(2.5)$ to those $M$ with $|M|=f$, i.e.,
$T_fP_{F_f}$, is perhaps more worth investigating. A typical example is the
additive divisor sum
$$
\sum_{n=1}^\infty\sigma_\lambda(n)\sigma_\mu(n+f)W(n/f),\eqno(2.7)
$$
where $\sigma_\lambda(n)=\sum_{d|n}d^\lambda$ with $\lambda\in\B{C}$; and
the smooth weight $W$ has a support in the positive reals. A detailed
discussion of the spectral decomposition of $(2.7)$ is developed in [27],
and the results there have recently found important applications in our 
joint work [19] (see also [20]) with M. Jutila, where is considered an instance of
the case $(5)$ above, and established a
new uniform bound for Hecke $L$-functions associated with cusp forms 
under $(0.2)$, which generalizes the classical subconvexity bound $(1.3)$ to a
vast family of
$L$-functions. 
\bigskip
\noindent 
{\bf 3.} In this and the next two sections, we shall deal with the Poincar\'e
series $P_F$, with an analogy of $(1.8)$ in mind. To begin with, we shall
collect elements of the theory of $\varGamma$-automorphic representations
of $G$, under the specification $(0.2)$:  Thus we write
$$
\r{n}[x]=\left[\matrix{1&x\cr&1}\right],\quad
\r{a}[y]=\left[\matrix{\sqrt{y}&\cr&1/\sqrt{y}}\right],\quad
\r{k}[\theta]=\left[\matrix{\phantom{-}\cos\theta&
\sin\theta\cr-\sin\theta &\cos\theta}\right].\eqno(3.1)
$$ 
Let $N=\left\{\r{n}[x]: x\in\B{R}\right\}$, $A=\left\{\r{a}[y]:  
y>0\right\}$, and
$K=\left\{\r{k}[\theta]:\theta\in\B{R}/\pi\B{Z}\right\}$, so that
$G=NAK$ be the Iwasawa decomposition of $G$.  We  
read it as $G\ni
\r{g}=\r{n}\r{a}\r{k}=\r{n}[x]\r{a}[y]\r{k}[\theta]$.
The Haar measures on the groups
$N$, $A$, $K$, $G$ are defined, respectively, by $d\r{n}=dx$,
$d\r{a}=dy/y$, $d\r{k}=d\theta/\pi$,
$d\r{g}=d\r{n}d\r{a}d\r{k}/y$,  with Lebesgue measures $dx$,
$dy$, $d\theta$. The space $L^2(\varGamma\backslash G)$ is composed of all
left $\varGamma$-automorphic functions on $G$, vectors for short, which
are square integrable over $\varGamma\backslash G$ against $d\r{g}$.
Elements of $G$ act unitarily on vectors from the right, and we have the
orthogonal decomposition into invariant subspaces
$$ 
L^2(\varGamma\backslash G)=\B{C}\cdot1\bigoplus
{}^c\!L^2(\varGamma\backslash G)\bigoplus
{}^e\!L^2(\varGamma\backslash G),\eqno(3.2)
$$ 
where ${}^c\!L^2$ is the cuspidal subspace, and
${}^e\!L^2$ is spanned by integrals of Eisenstein series.
\smallskip 
The cuspidal subspace splits into irreducible subspaces:
$$ 
{}^c\!L^2(\varGamma\backslash G)=\overline{\bigoplus V};\quad 
\Omega|_{V}=\left(\nu_V^2-\txt{1\over4}\right)\cdot1,
\eqno(3.3)
$$ 
where  
$\Omega=y^2\left(\partial_x^2+
\partial_y^2\right)-y\partial_x\partial_\theta$ is the Casimir
operator. Under
$(0.2)$, we can restrict our attention to two
cases: either $\nu_V\in i(0,\infty)$ or 
$\nu_V$ is equal to half a positive odd  
integer. According to the right action of $K$, the space
$V$ is decomposed into $K$-irreducible subspaces
$$ 
V=\overline{\bigoplus_{p=-\infty}^\infty V_p},\quad \dim V_p\le
1.\eqno(3.4)
$$ 
If it is not trivial, $V_p$ is spanned by a
$\varGamma$-automorphic function $\varphi_p$ such that
$\varphi_p(\r{g}\r{k}[\theta])=\exp(2ip\theta)\varphi(\r{g})$; it is called a
$\varGamma$-automorphic form of spectral parameter
$\nu_V$ and weight $2p$.
\smallskip 
Let us assume temporarily that $V$ belongs to the unitary
principal series, i.e., $\nu_V\in i(0,\infty)$ under $(0.2)$. 
Then $\dim V_p=1$ for all $p\in\B{Z}$, and there exists a complete
orthonormal system $\left\{\varphi_p\in V_p:\,p\in\B{Z}\right\}$ of $V$ 
such that
$$
\eqalignno{
\varphi_p(\r{g})&=\sum_{\scr{n=-\infty}\atop\scr{n\ne0}}
^\infty|n|^{-\nu}\varrho_V(n)
\e{A}_n\phi_p(\r{g};\nu_V)\cr
&=\sum_{\scr{n=-\infty}\atop\scr{n\ne0}}
^\infty{\varrho_V(n)\over\sqrt{|n|}}
\e{A}_{\sgn(n)}\phi_p(\r{a}[|n|]\r{g};\nu_V),
&(3.5)
}
$$ 
where $\phi_p(\r{g};\nu)=y^{{1/2}+\nu}\exp(2ip\theta)$, and
$\e{A}_u$  is the Jacquet operator:
$$
\e{A}_u\phi_p(\r{g};\nu)=\int_N
\exp(-2\pi iuv)\phi_p (\r{w}\r{n}[v]\r{g};\nu)d\r{n},\quad
\r{w}=\r{k}\left[\txt{1\over2}\pi\right].\eqno(3.6)
$$ 
It should be observed that the coefficients
$\varrho_V(n)$ in $(3.5)$ do not depend on the weight. We note that
for $u\in\B{R}^\times$
$$
\eqalignno{
\e{A}_u\phi_p(\r{g};\nu)&=y^{{1/2}-\nu}\exp(2\pi iux)\exp(2pi\theta)
\int_{-\infty}^\infty{\exp(2\pi yuv)
\over(v^2+1)^{{1\over2}+\nu}}
\left({v+i\over v-i}\right)^{ 
p}dv\cr &=(-1)^p\pi^{{1/2}+\nu}|u|^{\nu-1/2}
\exp(2\pi iux)\exp(2pi\theta) {W_{\sgn(u)p,\nu} (4\pi|u|
y)\over\Gamma(\sgn(u) p+{1\over2}+\nu)}, &(3.7)
}
$$ 
where $W_{\lambda,\mu}(y)$ is the Whittaker function (see [WW,
Chapter   XVI]). The first line is valid for
$\Re\nu>0$, while the second defines $\e{A}_n\phi_p$ for all
$\nu\in\B{C}$. In particular, we have the expansion
$$
\varphi_0(\r{g})={2\pi^{{1/2}+\nu_V}\over\Gamma({1\over2}+\nu_V)}
\sqrt{y}\sum_{\scr{n=-\infty}\atop\scr{n\ne0}}
^\infty\varrho_V(n)K_{\nu_V}(2\pi|n|y)
\exp(2\pi inx),\eqno(3.8)
$$ 
with $K_\nu$ being the $K$-Bessel function of order $\nu$. This is a
cusp-form on the hyperbolic upper half plane $G/K$.
\smallskip 
Next, let us consider a $V$ in the discrete series; that
is, $\nu_V=\ell-{1\over2}$, $1\le\ell\in\B{Z}$. We have, in place  
of $(3.4)$,
$$
\hbox{either\quad
$\displaystyle{V=\overline{\bigoplus_{p=\ell}^\infty   V_p}}$
\quad or\quad
$\displaystyle{V=\overline{\bigoplus_{p=-\infty}^{-\ell} V_p}}$}\,,
\eqno(3.9)
$$ 
with $\dim V_p=1$, corresponding to the holomorphic and the
antiholomorphic discrete series. The involution
$\r{g}=\r{nak}\mapsto\r{n}^{-1}\r{a}\r{k}^{-1}$ maps one to
the other.  In the holomorphic case, we have a complete orthonormal
system $\{\varphi_p:\, p\ge\ell\}$ in $V$ such that
$$
\varphi_V(\r{g})=\pi^{1/2-\ell}\left({\Gamma(p+\ell)\over
\Gamma(p-\ell+1)}\right)^{1/2}\sum_{n=1}^\infty n^{-\nu_V}\varrho_V(n)
\e{A}_n\phi_p(\r{g};\nu_V).\eqno(3.10)
$$
In particular, we have
$$
\varphi_\ell(\r{g})=(-1)^\ell{2^{2\ell}\pi^{{1/2}+\ell}\over
\sqrt{\Gamma(2\ell)}}\exp(2i\ell\theta)y^\ell
\sum_{n=1}^\infty\varrho_V(n)n^{\ell-{1/2}}\exp(2\pi in(x+iy)).
\eqno(3.11)
$$ 
This infinite sum is a holomorphic cusp-form of weight $2\ell$ on $G/K$. 
\medskip 
With this, we may assume further that all $V$ be Hecke invariant, so that
there exists, for any integer $n\ge1$, a real number $\tau_V(n)$ such that
$$ 
T_n|_V=\tau_V(n)\cdot1,\eqno(3.12)
$$ 
where $T_n$ is as in $(2.6)$. Thus, for any non-zero integer $n$,
$$
\varrho_V(n)=\varrho_V(\sgn(n))\tau_V(|n|).\eqno(3.13)
$$ 
We may introduce the convention that $\varrho_V(-1)=0$ and
$\varrho_V(1)=0$ for $V$ in the holomorphic and antiholomorphic
discrete series, respectively, and $\varrho_V(-1)=\epsilon_V
\varrho_V(1)$ with $\epsilon_V=\pm1$ for $V$ in the unitary principal
series. We associate, with each $V$, the Hecke series
$$ 
H_V(s)=\sum_{n=1}^\infty \tau_V(n)n^{-s},\quad \Re s>1,\eqno(3.14)
$$ 
which continues to an entire function.
\bigskip
\noindent
{\bf 4.} With this, we are going to decompose $P_F$ via $(3.2)$. We may
restrict ourselves to the orthogonal projection
$\varpi_V$ with $V$  in the unitary principal series. In fact, 
the discrete series is highly analogous,  the projection to the subspace 
${}^e\!L^2(\varGamma\backslash G)$ is facilitated by  that
Eisenstein series are explicitly defined, and the space of
constant functions gives no specific problem.
\smallskip
We are of course concerned with how to compute $\varpi_VP_F$ in
terms of $F$, or more generally with the harmonic analysis on $V$ 
while the variable $\r{g}\in G$ is restricted to the
big Bruhat cell. To this end we shall employ a reasoning 
which we term the Kirillov scheme. This is because we utilize an operator
$\e{K}$, defined by $(4.3)$ below, that apparently originated 
in A.A. Kirillov [21].  We shall exploit two basic properties of $\e{K}$, 
and they are embodied here in two lemmas, respectively:
\medskip
\noindent
{\bf Lemma 1.} {\it Let $U=U_\nu$ with $\nu\in i\B{R}$ be the Hilbert space
$$
\overline{\bigoplus_{p=-\infty}^\infty\B{C}\phi_p},\quad
\phi_p(\r{g})=\phi_p(\r{g};\nu),\eqno(4.1)
$$
equipped with the norm
$$
\Vert{\phi}\Vert_{U}=\sqrt{\sum_{p=-\infty}^\infty|c_p|^2},\quad
\phi=\sum_{p=-\infty}^\infty c_p\phi_p.\eqno(4.2)
$$
For $u\in\B{R}^\times$ and smooth $\phi\in U$, that is, with $c_p$
decaying faster than any negative power of $|p|$, we let
$$
\e{K}\phi(u)=|u|^{1/2-\nu}\e{A}_u\phi(1)=\e{A}_{\sgn(u)}\phi(\r{a}[|u|]).
\eqno(4.3)
$$ 
Then the operator $\e{K}$ maps $U$ unitarily onto
$L^2(\B{R}^\times,(1/\pi)d^\times)$, where $d^\times\! u=du/|u|$. 
}
\medskip
\noindent{\bf Lemma 2.} {\it Let us define the Bessel function of
representations of $\r{PSL}_2(\B{R})$ as to be
$$
j_\nu(u)=\pi{\sqrt{|u|}\over\sin\pi\nu}\left(J_{-2\nu}^{\sgn(u)}
(4\pi\sqrt{|u|})-J^{\sgn(u)}_{2\nu}(4\pi\sqrt{|u|})\right),
\eqno(4.4)
$$ 
with $J^+_\nu=J_\nu$ and $J^-_\nu=I_\nu$ in the ordinary 
notation for Bessel functions. Then, for any smooth $\phi\in U_\nu$, we have
$$
\e{K}R_\r{w}\phi(u)=\int_{\B{R}^\times}j_\nu(uv)\e{K}
\phi(v)d^\times\!v,\quad u\in\B{R}^\times.\eqno(4.5)
$$
}
\medskip
\noindent
{\it Proof\/}. It should be stressed that the definition $(4.3)$ is taken from
[9][35], and somewhat different from that employed in [32][33].
A proof of the unitarity of $\e{K}$ is given in [32] [33, Theorem 1]. It depends
on the following  integral formula via the second line of $(3.7)$: For any
$\alpha,\beta\in\B{C}$ and $|\Re\nu|<{1\over2}$
$$
\eqalignno{
\qquad\int_0^\infty& W_{\lambda,\nu}(u) W_{\mu,\nu}(u){du\over u}=
{\pi\over(\lambda-\mu)\sin(2\pi\nu)}\cr
&\times\left[{1\over\Gamma({1\over2}-\lambda+\nu)\Gamma({1\over2}
-\mu-\nu)}-{1\over\Gamma({1\over2}-\lambda-\nu)\Gamma({1\over2}
-\mu+\nu)}\right], &(4.6)
}
$$ 
which is tabulated as [14, eq.\ 7.611(3)]. The verification of this made in [32]
[33]  employs the Whittaker differential equation ([41, p.\ 337]) which is related to
the Casimir operator. The surjectivity of $\e{K}$ is proved in [9] via  the first
line of $(3.7)$ and the completeness of the system
$\left\{((v+i)/(v-i))^p:\,p\in\B{Z}\right\}$ in the space
$L^2\left(\B{R},dv/(\pi(v^2+1))\right)$.  As to Lemma 2, the
realization $(4.5)$ of the action of $\r{w}$ the Weyl element in terms of  the
space $L^2(\B{R}^\times,(1/\pi)d^\times)$ seems to have been published for the
first time by N.Ja.\ Vilenkin (see [38, Section 7 of Chapter VII] as well as [39,
eq.\ $(17)$ on p.\ 454]), though the concept of the Bessel function of
representations had been coined by I.M. Gel'fand, M.I. Graev and I.I.
Pyatetski-Shapiro [12]. Two independent proofs are known; they  are
conceptually different. One is due to M. Baruch and Z. Mao [2], which is along
the line of [38] and fills a gap therein concerning a convergence
issue. The other is due to ourselves [31][32][33], and seems more in line with the
purpose of the present work. It is shown there that $(4.5)$ is in fact
equivalent to the Jacquet--Langlands local functional equation ([17,
Theorem 5.15])
$$
\eqalignno{ 
(-1)^p\Gamma_p(s) =&2^{1-2s}\pi^{-2s}\Gamma(s+\nu)
\Gamma(s-\nu)\cr &\times\left(\cos(\pi s)\Gamma_p(1-s)+\cos(\pi\nu)
\Gamma_{-p}(1-s)\right)&(4.7)  
}
$$ 
for the Mellin transform
$$
\Gamma_p(s)=\int_0^\infty
\e{A}_u\phi_p(1)u^{s-\nu-1}du,\eqno(4.8)
$$ 
which continues meromorphically to $\B{C}$. The Mellin inversion of $(4.7)$
coupled with $(4.9)$--$(4.10)$ below gives  $(4.5)$ for $\phi=\phi_p$; the
extension to any smooth $\phi\in U$ is easy. Thus, if $|\Re\nu|-{1\over2}<\Re
s$, then
$$
\int_{-\infty}^0 j_\nu(u)|u|^{s-1}du={1\over\pi}(2\pi)^{-2s}\cos(\pi
\nu)\Gamma\left(s+\txt{1\over2}+\nu\right)
\Gamma\left(s+\txt{1\over2}-\nu\right);\eqno(4.9)
$$ 
and if $|\Re\nu|-{1\over2}<\Re s<-{1\over4}$, then
$$
\int_0^\infty j_\nu(u)u^{s-1}du=-{1\over\pi}(2\pi)^{-2s}\sin(\pi s)
\Gamma\left(s+\txt{1\over2}+\nu\right)
\Gamma\left(s+\txt{1\over2}-\nu\right);\eqno(4.10)
$$ 
The former follows from [40, eq.\ $(8)$ in Section 13.21], and the
latter from [ibid, eq.\ $(1)$ in Section 13.24]. This ends our brief
discussion on the proof of the above lemmas. 
\smallskip
The above extends not only to the discrete
and the complementary series  but also to the complex situation,
i.e., to $\r{PSL}_2(\B{C})$, as is to be shown in Section 7.
\bigskip
\noindent
{\bf 5.} Now we shall carry out the computation of $\varpi_VP_F$ via the
Kirillov scheme. It should be stressed that absolute convergence
required below can readily be confirmed, provided $F$ is sufficiently
smooth.
\smallskip
Thus, the projection to $V_p$ is, by the unfolding argument,
$$
\eqalignno{\langle P_F,
\varphi_p\rangle_{\varGamma\backslash G}&=\int_G
F(\r{g})\overline{\varphi_p(\r{g})}d\r{g}\cr
&=\overline{\varrho_V(1)}
\sum_{m=1}^\infty{\tau_V(m)\over\sqrt{m}}
(\Phi_p^++\epsilon_V\Phi_p^-)F_m(\nu_V),&(5.1) }
$$ 
where $(3.5)$, $(3.13)$ are used; $F_m(\r{g})=F(\r{a}[m]^{-1}\r{g})$ and
$$
\Phi_p^\delta F(\nu)=\int_G
F(\r{g})\overline{\e{A}_\delta\phi_p(\r{g})}d\r{g}.\eqno(5.2)
$$ 
Thus
$$
\eqalignno{
\varpi_V P_F(\r{g})&=\sum_{p=-\infty}^\infty\langle P_F,
\varphi_p\rangle_{\varGamma\backslash G}\varphi_p(\r{g})\cr
&=|\varrho_V(1)|^2\sum_{m=1}^\infty
\sum_{n=1}^\infty{\tau_V(m)\tau_V(n)\over\sqrt{mn}}\cr 
&\times\left(\e{B}^{(+,+)}+\e{B}^{(-,-)}+\epsilon_V\e{B}^{(+,-)}
+\epsilon_V\e{B}^{(-,+)}\right)F_m(\r{a}[n]\r{g};\,\nu_V),&(5.3) 
}
$$ 
where
$$
\eqalignno{
\e{B}^{(\delta_1,\delta_2)}
F(\r{g};\nu)&=\sum_{p=-\infty}^\infty\Phi_p^{\delta_1}F(\nu)\,
\e{A}_{\delta_2}\phi_p(\r{g};\nu)\cr &=\exp(2\pi i\delta_2
x)\sum_{p=-\infty}^\infty\Phi_p^{\delta_1} 
F(\nu)\e{A}_{\delta_2}\phi_p(\r{a}[y])\exp(2ip\theta).&(5.4) 
}
$$ 
\par
Since it can be asserted, with an appropriate change of $F$, 
that our interest is in the value
$\varpi_V P_F(1)$ (see $(2.5)$), we may restrict ourselves to the subgroup
$A$. Namely, it suffices to consider
$\e{B}^{(\delta_1,\delta_2)}F(\r{a}[y];\nu)$, with a new $F$; and this  can be
expressed in terms of the Kirillov  operator: 
$$
\e{B}^{(\delta_1,\delta_2)}F(\r{a}[y];\nu)=\e{K}\e{L}^{\delta_1}F(\delta_2y),
\quad \e{L}^\delta F=\sum_{p=-\infty}^\infty\Phi_p^\delta F(\nu)\phi_p,
\eqno(5.5)
$$
where $\e{L}^\delta F\in U$ is smooth.
Hence, we may proceed in the sense of weak convergence, 
appealing to Lemma 1; and we observe
$$
\Phi_p^{\delta}F(\nu)=\langle \e{L}^\delta
F,\phi_p\rangle_{U} ={1\over\pi}\int_{\B{R}^\times}
\e{K}\e{L}^\delta F(u)\overline{\e{K}\phi_p(u)} d^\times\!
u.\eqno(5.6)
$$ 
This means that if we are able to transform $(5.2)$ into
$$
\Phi_p^{\delta}F(\nu)={1\over\pi}\int_{\B{R}^\times}Y^\delta(u)
\overline{\e{K}\phi_p(u)}d^\times\! u,\eqno(5.7)
$$ 
then it should follow that
$$
\e{B}^{(\delta_1,\delta_2)}F(\r{a}[y];\nu)=Y^{\delta_1}(\delta_2
y),\eqno(5.8)
$$ 
because of the surjectivity assertion in the lemma. 
\smallskip
Since the integral in $(5.2)$ is in fact over the big Bruhat cell, we
perform the change of variables accordingly. We have instead
$$
\Phi_p^{\delta}F(\nu) =\int_0^\infty\int_{N\r{w}N}F(\r{a}[u]\r{g})
\overline{R_\r{g}\e{A}_\delta\phi_p(\r{a}[u])}d\dot\r{g}{du\over
u}.\eqno(5.9)
$$ 
Here $R_\r{g}$ is the right translation with $\r{g}=\r{n}[x_1]\r{w}\r{n}[x_2]$,
and $d\dot\r{g}=dx_1dx_2/\pi$.  We observe that
$$
\eqalignno{ 
R_\r{g}\e{A}_\delta\phi_p(\r{a}[u])&=\exp(2\pi i\delta
x_1u)\e{A}_\delta
R_\r{w}R_{\r{n}[x_2]}\phi_p(\r{a}[u])\cr
&=\exp(2\pi i\delta x_1u)\e{K}R_\r{w}R_{\r{n}[x_2]}\phi_p(\delta
u).&(5.10)
}
$$ 
By Lemma 2 this is replaced by
$$
R_\r{g}\e{A}_\delta\phi_p(\r{a}[u])=
\exp(2\pi i\delta x_1u)\int_{\B{R}^\times} \exp(2\pi
ix_2v)j_\nu(\delta uv)\e{K}\phi_p(v)
d^\times\!v,\eqno(5.11)
$$ 
and $(5.9)$ by
$$
\eqalignno{
\Phi_p^{\delta}F(\nu)=&{1\over\pi}\int_0^\infty\int_{\B{R}^2}
F(\r{a}[u]\r{n}[x_1]\r{w}\r{n}[x_2])\exp(-2\pi i\delta x_1u)\cr
&\times\int_{\B{R}^\times} \exp(-2\pi ix_2v)j_\nu(\delta
uv)\overline{\e{K}\phi_p(v)}d^\times\!v 
dx_1dx_2{du\over u}.&(5.12) }
$$ 
Hence we find via $(5.8)$ that
$$
\eqalignno{ &\e{B}^{(\delta_1,\delta_2)}F(\r{a}[y];\nu)
=\int_0^\infty j_\nu(\delta_1\delta_2 yu)\cr &\times\left\{ \int_{\B{R}^2}
F(\r{a}[u]\r{n}[x_1]\r{w}\r{n}[x_2])\exp(-2\pi i\delta_1u x_1-2\pi
i\delta_2yx_2)dx_1dx_2\right\}{du\over u},&(5.13)  
}
$$ 
which ends the application of the Kirillov scheme. 
\medskip
We may compare $(4.7)$ with $(1.7)$. Then $(4.5)$ may 
also be compared with $(1.8)$. That is, the formula $(5.3)$ coupled with
$(5.13)$ which is a local assertion derived from $(4.5)$ corresponds to
$(1.8)$. As remarked above already, analogues of Lemmas 1 and 2 are
shown in [9][32] for the discrete series representations and the complementary
series, although the latter is irrelevant under $(0.2)$. Hence,
returning to $(3.2)$, we obtain a genuine extension of $(1.8)$. The final result
is, however, too complicated to be stated as an independent assertion. We
should instead be content with the local expression $(5.13)$ and with the fact
that we have found that the combination of Lemmas 1 and 2 is the key
implement.
\smallskip
Hence, what corresponds to the cosine-transform in $(1.8)$ is $(5.3)$ 
with $(5.13)$. One may desire to compute the double sum $(5.3)$ and the last
double integral into closed forms. In the applications to $\e{M}(\zeta^2,g)$ and
to the sum
$(2.7)$, which will be briefly dealt with below, we are in a fortuitous situation
that the double sum is transformed into a product of two values of $H_V$. As
to the double integral, it is a Fourier transform over the Euclidean plane, and
thus might be expressed in terms of a Bessel transform. With
$\e{M}(\zeta^2,g)$ as well as $(2.7)$, the situation turns out in fact to be as
such, and we shall see that $(5.13)$ is expressed as an integral transform
whose kernel is a convolution of two instances of the Bessel function of
representations (see
$(6.2)$ below). 
\smallskip
Thus the matter  seems to depend much on the specific
nature of the seed $F$.  Nevertheless, with any smooth $F$, one
might appeal to Mellin transform of several variables, and 
the above could be pushed into a more closed form. 
\bigskip
\noindent
{\bf 6.} We are now at the stage to render the spectral decomposition of
$\e{M}(\zeta^2,g)$ in terms of notions from representation theory:  Thus, let us
put
$$
\eqalignno{
\Theta(\nu;g)&={1\over4\cos(\pi\nu)}\int_0^\infty
\left({u\over u+1}\right)^{1/2}g_c\left(
\log\left(1+{1/u}\right)\right)
\Xi(u;\nu)d^\times\! u,&(6.1)\cr
\Xi(u;\nu)&=\int_{\B{R}^\times} j_0(-v)j_\nu\left({v\over u}\right)
{d^\times\! v\over\sqrt{|v|}}.&(6.2)
}
$$
Then we have
$$
\e{M}(\zeta^2,g)=\left\{\e{M}^{(r)}+\e{M}^{(c)}
+\e{M}^{(e)}\right\}(\zeta^2,g),\eqno(6.3)
$$ 
where
$$
\eqalignno{
\e{M}^{(c)}(\zeta^2,g)&=\sum_{V}
\alpha_VH_V\left(\txt{1\over2}\right)^3
\Theta(\nu_V;g),&(6.4)\cr
\e{M}^{(e)}(\zeta^2,g)&=\int_{(0)}
{\left|\zeta\left({1\over2}+\nu\right)\right|^6
\over|\zeta(1+2\nu)|^2}\Theta(\nu;g){d\nu\over2\pi
i},&(6.5) 
}
$$ 
with $\alpha_V=|\varrho_V(1)|^2+|\varrho_V(- 
1)|^2$. The $V$ runs over a maximal orthogonal system of
Hecke-invariant cuspidal $\varGamma$-automorphic representations of
$G$. Apart the term
$2\pi\Re\{(\log(2\pi)-\gamma_E)g({1\over2}i)-
{1\over2}ig'({1\over2}i)\}$,  the $\e{M}^{(r)}(\zeta^2,g)$ is an
integral transform of $g$ whose kernel is given explicitly in terms of
logarithmic derivatives of the Gamma function.
\smallskip
The proof of $(6.3)$ as developed in [9] starts with the integration of
$\zeta(z_1+it)\zeta(z_2+it)\zeta(z_3-it)\zeta(z_4-it)$ against $g(t)dt$ over
$\B{R}$, where $(z_1,z_2,z_3,z_4)$ is to remain in the region of absolute
convergence. The device $(2.3)$--$(2.5)$ is to be applied. However, the naive
choice of the seed does not work well, because we require it be
smooth and of rapid decay on $G$. Such a choice of the seed
is not difficult but somewhat subtle. Thus, we are forced to employ instead a
sequence of suitable $F$'s and a limiting procedure with respect to $F$.
Nevertheless, those $F$ chosen in [9] are nice in that they are
$A$-equivariant, i.e., $F(\r{a}[y]\r{g})=y^\omega F(\r{g})$ with an
$\omega\in\B{C}$ being independent of $F$. Hence the summation over $f$ in
$(2.5)$ is the same as the multiplication by a value of
$H_V$ at each $V$; and also the sum  over $m$ in $(5.1)$ can be written in
terms of a product of a value of $H_V$ and
$(\Phi_p^++\epsilon_V\Phi_p^-)F(\nu_V)$. In particular, the omission  of $m$
in $(5.4)$ is possible without changing the specification of
$F$, which amounts to a considerable simplification in the subsequent
discussion leading to $(5.13)$. Besides, this makes $(5.13)$ easier to handle;
taking the limit in $F$ we come already close to the expression $(6.2)$.
Moreover, the sum over $n$ in $(5.3)$ yields now another factor in a value of
$H_V$. Then it remains to  perform analytic continuation and specialization
with respect to $(z_1,z_2,z_3,z_4)$.
Thus the factor $j_0$ in $(6.2)$ stands for the integral kernel of the Bessel
transform that emerges from the inner double integral of $(5.13)$ at the end of
the whole procedure (see [9, (7.34)]). This explains how $(6.2)$ originates and
reveals especially the mechanism behind $(6.4)$.
\medskip
The same argument can be applied to the additive divisor sum $(2.7)$. We are
led to a spectral decomposition analogous to $(6.3)$. Specifically, the cuspidal
part is found to be
$$
\eqalignno{
{1\over4}f^{(\lambda+\mu+1)/2}
\sum_V\alpha_V&\tau_V(f)H_V\left(\txt{1\over2}(1-\lambda-\mu)\right)
H_V\left(\txt{1\over2}(1+\lambda-\mu)\right)\cr
&\times\left(\Psi_++\epsilon_V\Psi_-\right)
(\nu_V;\lambda,\mu;W), &(6.6)
}
$$
where
$$
\Psi_\delta(\nu;\lambda,\mu;W)=\int_0^\infty W(u)
\Lambda_\delta(u;\nu;\lambda,\mu)u^{(\lambda+\mu)/2+1}d^\times\!
u,\eqno(6.7)
$$
with
$$
\Lambda_\delta(u;\nu;\lambda,\mu)=\int_0^\infty j_{\lambda/2}(-\delta v)
j_\nu(\delta v/u){d^\times\! v\over v^{(\mu+1)/2}}.\eqno(6.8)
$$
It is safe to keep both $|\Re\lambda|,\,|\Re\mu|$
sufficiently small so that $(6.6)$ holds with the
expression $(6.8)$. However, as can be seen from 
$(4.9)$--$(4.10)$, $\Lambda_\delta$ can be expressed in terms of
the Mellin inversion of a product of four Gamma factors, and then 
$(6.7)$ allows us to continue $(6.6)$ analytically to quite a wide
domain of $(\lambda,\mu)$.
\medskip
It has been explained in the above how the factor $j_0$ in
$(6.2)$ turns up; the same can be applied to $j_{\lambda/2}$ in $(6.8)$.
However, it is not done in any framework of metric theory. This  is
sharply different from the situation with another factor $j_\nu$ shared
in both equations. Thus it remains still to find a genuine characterization of
the factors $j_0$,
$j_{\lambda/2}$. In passing, we note that instead of $(6.2)$ we may write
$$
\Xi(u;\nu)=2\Re\Bigg\{u^{-1/2-\nu}\left(1-{1\over\sin(\pi\nu)}\right)
{\Gamma^2({1\over2}+\nu)\over\Gamma(1+2\nu)}{}_2F_1
\left(\txt{1\over2}+\nu,\txt{1\over2}+\nu;
1+2\nu;-1/u\right)\Bigg\},\eqno(6.9)
$$
with the Gaussian hypergeometric function ${}_2F_1$ (see [30, (4.7.2)]). 
This reminds us of the free-space resolvent kernel of the hyperbolic Laplacian (see
[ibid, (1.1.49)]), a fact that appears mysterious to us.
\medskip
It might be expedient to make here a digression on a historical background:
A prototype of the spectral decomposition of $\e{M}(\zeta^2,g)$ was obtained
by the present author in [24][25], which was afterwards improved to $(6.3)$ in
[30, Theorem 4.1]. However, the assertion there did not reach the expression
$(6.2)$; it was stated with $(6.9)$. A reason for this is in that there we used
the Kloosterman-Spectral sum formula of N.V. Kuznetsov [30, Theorems 2.3 and
2.5], which is pretty handy but hides the mechanism working behind the integral
transform appearing on the spectral side. Note that in the above we 
dispensed with  Kuznetsov's sum formula. The argument of [9], whose
most salient part is depicted in the previous section, is admittedly more
involved than that in [30], but this is much due to the fact that we started
from the very fundamental assertion $(3.2)$,
whereas the discussion in [30] lacks the perspective offered by representation
theory. Thus, the appearance in
Kuznetsov's sum formula and consequently in [30, Theorem 4.1] of the
contribution of holomorphic cusp forms was just an accidental  byproduct of a
technical marvel and remained mysterious there. Our discussion of
$\e{M}(\zeta^2,g)$ in terms of a Poincar\'e series on the group $G$ allows us
to see all contributions of cusp forms in a fairly equal term, since our method
is based on $(3.2)$, where all irreducible representations have equal rights.
\smallskip
It was Bruggeman [3][4] who tried for the first time to understand, via
$(3.2)$, all the terms on the spectral side in Kuznetsov's sum formula. However,
the real comprehension of the structure supporting the sum formula appears to have
been done by J.W. Cogdell and I. Piatetski-Shapiro in [11]. In
particular, the Kirillov scheme together with the r\^ole of the Bessel function
of representations was developed there, and Kuznetosv's sum
formula was newly proved, though their discussion appears 
sketchy to us.  The authors of [9] were inspired by
the work [11].
\bigskip
\noindent
{\bf 7.} The aim of this and the next sections is to show that the above
discussion extends to the situation $(0.3)$. In particular, we are going to show
the complex analogues of Lemmas 1 and 2. Note that some symbols used under
$(0.2)$ are now assigned to corresponding notions under $(0.3)$; this
convention should not cause any confusion.
\smallskip
Thus, let $G=\r{PSL}_2(\B{C})$, and put
$$
\r{n}[z]=\left[\matrix{1&z\cr&1}\right],\;\r{h}[u]=\left[
\matrix{u&\cr&1/u}\right],\;\r{k}=
\left[\matrix{\alpha&\beta\cr-\overline{\beta}&
\overline{\alpha}}\right],\eqno(7.1)
$$ 
where $z,u,\alpha,\beta\in\B{C}$ with $u\ne0$,
$|\alpha|^2+|\beta|^2=1$; and also
$N=\{\r{n}[z]:z\in\B{C}\}$, $A=\{\r{a}[r]:r>0\}$, 
$K=\r{PSU}(2)=\{\r{k}[\alpha,\beta]:\alpha,\beta\in\B{C}\}$ 
with $\r{a}[r]=\r{h}[\sqrt{r}]$. In terms of the Euler angles
$\varphi,\theta,\psi$, we have
$$
\r{k}[\alpha,\beta]=\r{h}[e^{i\varphi/2}]\r{v}[i\theta]\r{h}[e^{i\psi/2}],\quad
\r{v}[\theta]=\r{k}[\cosh(\theta/2),
\sinh(\theta/2)].\eqno(7.2)
$$ 
The Iwasawa decomposition $G=NAK$ is read as $G\ni\r{g}=\r{n}[z]
\r{a}[r]\r{k}[\alpha,\beta]$. The Haar measures on respective
groups are given by $d\r{n}=dz$, $d\r{a}=dr/r$, $d\r{k}=
\sin\theta d\varphi d\theta d\psi/(8\pi^2)$, and
$d\r{g}=d\r{n}d\r{a}d\r{k}/r^2$.  With this, the Hilbert space
$L^2(\varGamma\backslash G)$ is formed, to which $G$ acts from
the right; and we have an exact analogue of $(3.2)$. See [5][8]
for more details of what follows.  
\smallskip
The cuspidal
subspace decomposes into irreducible subspaces
$$ 
{}^c\!L^2(\varGamma\backslash G)=\overline{\bigoplus
V}.\eqno(7.3)
$$
To classify representations $V$, we need two Casimir operators
$\Omega_\pm$, $\Omega_-=\overline{\Omega_+}$, where
$$
\eqalignno{
\Omega_+={1\over2}r^2&\partial_z\partial_{\overline{z}}+
{1\over2}re^{i\varphi}\cot\theta\partial_z\partial_\varphi
-{1\over2}ire^{i\varphi}\partial_z\partial_\theta-
{re^{i\varphi}\over2\sin\theta}\partial_z\partial_\psi
\cr&+{1\over8}r^2\partial_r^2-{1\over4}ir
\partial_r\partial_\varphi
-{1\over8}\partial_\varphi^2
-{1\over8}r\partial_r+{1\over4}
i\partial_\varphi. &(7.4)
}
$$ 
They become constant multiplications in each $V$:
$$
\Omega_\pm\vert_V=\chi_V^\pm\cdot1=
{1\over8}((p_V\mp\nu_V)^2-1)\cdot1,\quad p_V\in\B{Z},\;\nu_V\in
i[0,\infty).\eqno(7.5)
$$ 
The pair $(p_V,\nu_V)$ is called the spectral parameter of $V$, but in
the sequel we shall write simply $(p_V,\nu_V)=(p,\nu)$.
\smallskip
According to the action of $K$, the space $V$
decomposes into $K$-irreducible subspaces
$$ 
V=\overline{\bigoplus_{|p|\le l,|q|\le l} V_{l,q}},\quad \dim
V_{l,q}=1.\eqno(7.6)
$$ 
To describe this precisely, let $\Omega_K$ be the Casimir
element of the universal enveloping algebra of $K$ defined by
$$
\Omega_K={1\over2\sin^2\theta}\left(\partial^2_\varphi
+\sin^2\theta\partial^2_\theta+\partial^2_\psi
-2\cos\theta\partial_\varphi\partial_\psi+\sin\theta\cos\theta\,
\partial_\theta\right).\eqno(7.7)
$$ 
Then 
$$ 
V_{l,q}=\left\{F\in V:
\Omega_K F=-\txt{1\over2}l(l+1),\,\partial_\psi F
=-iqF\right\}.\eqno(7.8)
$$  
Any non-zero element of $V_{l,q}$ is called a
$\varGamma$-automorphic form of spectral parameter $(p,\nu)$
and $K$-type $(l,q)$. 
\smallskip 
Next, we define functions $\Phi^l_{p,q}$ on $K$ by
$$ 
(\alpha X-\overline{\beta})^{l-q}(\beta
X+\overline{\alpha})^{l+q}
=\sum_{p=-l}^l\Phi_{p,q}^l(\r{k}[\alpha,\beta])X^{l-p}.
\eqno(7.9)
$$
The system $\{\Phi_{p,q}^l: |p|,|q|\le l, 1\le l\}$ is a complete
orthogonal basis of $L^2(K)$ with norms
$$
\Vert\Phi_{p,q}^l\Vert_K={1\over\sqrt{l+{1\over2}}}
{2l\choose l-p}^{1/2}{2l\choose l-q}^{-1/2}.\eqno(7.10)
$$ 
Here are some of its properties which we shall need later: Under the convention
that $\Phi_{p,q}^l\equiv0$ if the condition $|p|,|q|\le l$ is violated, we have
$$
\eqalignno{
\Phi_{p,q}^l(\r{k}[\alpha,\beta])&=e^{-ip\varphi-iq\psi}
\Phi_{p,q}^l(\r{v}[i\theta]),&(7.11)\cr
2\partial_\theta\Phi_{p,q}^l(\r{v}[i\theta])&=
i(l+p+1)\Phi_{p+1,q}^l(\r{v}[i\theta])
+i(l-p+1)\Phi_{p-1,q}^l(\r{v}[i\theta])\cr
&=i(l-q)\Phi_{p,q+1}^l(\r{v}[i\theta])
+i(l+q)\Phi_{p,q-1}^l(\r{v}[i\theta]), &(7.12)\cr
2{p-q\cos\theta\over\sin\theta}\Phi_{p,q}^l(\r{v}[i\theta])
&=i(l-q)\Phi_{p,q+1}-i(l+q)\Phi_{p,q-1}^l(\r{v}[i\theta]),&(7.13)\cr
2{q-p\cos\theta\over\sin\theta}\Phi_{p,q}^l(\r{v}[i\theta])
&=i(l+p+1)\Phi_{p+1,q}^l-i(l-p+1)\Phi_{p-1,q}^l(\r{v}[i\theta]).&(7.14)
}
$$
For a verification of $(7.12)$--$(7.14)$ see [5, Lemma 5].
\smallskip
With this, we put $\phi_{l,q}(\r{g};\nu)={r^{1+\nu}\Phi^l_{p,q}(\r{k})
/\Vert\Phi_{p,q}^l\Vert_K}$; and its Jacquet transform is defined by
$$
\e{A}_u\phi_{l,q}(\r{g};\nu)=\int_N
\exp(-2\pi i\Re(uv))
\phi_{l,q}(\r{w}\r{n}[v]\r{g};\nu)d\r{n}.\eqno(7.15)
$$
Let $\varphi_{l,q}$ be a generating vector of
$V_{l,q}$, so that $\{\varphi_{l,q}:|p|\le l, |q|\le l\}$ is a
complete orthonormal system in $V$. Then we have
$$
\varphi_{l,q}(\r{g})=\sum_{\omega\ne0}
|\omega|^{-\nu}\varrho_V(\omega)\e{A}_\omega\phi_{l,q}(\r{g};\nu),\quad
\omega\in\B{Z}[i],\eqno(7.16)
$$
which is precisely an analogue of the first line of $(3.5)$. 
\smallskip 
Further, we define the Hecke operator for each non-zero $f\in\B{Z}[i]$ by
$$
T_fF(\r{g})={1\over4|f|}
\sum_{d|f}\;\sum_{b\bmod\,d}F\left(
\left[\matrix{\sqrt{f}/d&b\cr&d/\sqrt{f}}\right]
\r{g}\right);\eqno(7.17)
$$
and we assume that all $V$ are Hecke invariant so that there exists a
real number $\tau_V(f)$ such that $T_f|V=\tau_V(f)\cdot 1$.
In particular we have, for all non-zero $n\in\B{Z}[i]$,
$$
\varrho_V(n)=\varrho_V(1)(n/|n|)^p\tau_V(n),\eqno(7.18)
$$
We have
$$
 \tau_V(-n)=\tau_V(n),\quad
\tau_V(in)=\epsilon_V\tau_V(n),\quad\epsilon_V=\pm1\eqno(7.19)
$$
The Hecke $L$-function of the space $V$ is defined by
$$
H_V(s)={1\over4}\sum_{n\ne0}\tau_V(n) |n|^{-2s},\eqno(7.20)
$$
which continues to an entire function; note that $H_V\equiv0$ whenever
$\epsilon_V=-1$.
\medskip
Now, the complex analogues of Lemmas 1 and 2 are as follows:
\medskip
\noindent
{\bf Lemma 3.} {\it  Let $U=U_{p,\nu}$ be the Hilbert space
$$
\overline{\bigoplus_{|p|\le l, |q|\le l}\B{C}\phi_{l,q}},
\quad \phi_{l,q}(\r{g})=\phi_{l,q}(\r{g};\nu)\eqno(7.21)
$$
equipped with the norm
$$
\Vert{\phi}\Vert_{U}=\sqrt{\sum_{|p|\le l, |q|\le l}|c_{l,q}|^2},\quad
\phi=\sum_{|p|\le l, |q|\le l} c_{l,q}\phi_{l,q}\,.\eqno(7.22)
$$
For $u\in\B{C}^\times$ and smooth $\phi\in U$, we let
$$
\e{K}\phi(u)=|u|^{1-\nu}(u/|u|)^p\e{A}_u\phi(1).\eqno(7.23)
$$
Then the operator $\e{K}$ maps $U$ unitarily onto
$L^2(\B{C}^\times,(2/\pi)d^\times)$, where $d^\times\! u= du/|u|^2$. 
}
\medskip
\noindent 
{\bf Lemma 4.}  {\it 
Let us define the Bessel function of representations of $\r{PSL}_2(\B{C})$
as to be
$$ 
j_{p,\nu}(u)=2\pi^2{|u|^2\over\sin(\pi\nu)}\left(
J_{p-\nu}(2\pi u)J_{-p-\nu}(2\pi\overline{u})-
J_{-p+\nu}(2\pi u)J_{p+\nu}(2\pi\overline{u})\right).\eqno(7.24)
$$ 
Then we have, for any smooth $\phi\in U_{p,\nu}$, 
$$ 
\e{K}R_\r{w}\phi(u^2)=\int_{\B{C}^\times}j_{p,\nu}(uv)\e{K}\phi(v^2)
d^\times\! v,\quad u\in\B{C}^\times.\eqno(7.25)
$$ 
Here $J_{p-\nu}(u)J_{-p-\nu}(\overline{u})$ is understood to be equal to 
$(u/|u|)^p|u|^{-2\nu}J^*_{p-\nu}(u)J^*_{-p-\nu}(\overline{u})$ where $J_\mu^*(u)$
is the entire function that coincides with $u^{-\mu}J_\mu(u)$ when $u>0$.
}
\medskip
With this, the analogues of $(6.1)$--$(6.5)$ for the Dedekind zeta-function
$\zeta_k$ of the Gaussian number field $k=\B{Q}(i)$ can be rendered as
follows: Let us put
$$
\eqalignno{
\Theta(p,\nu;g)&={\nu\over16\sin(\pi\nu)}\int_\B{C}{|u|\over
|u+1|}
g_c\left(2\log|1+1/u|\right)\Xi(u;p,\nu)d^\times\! u,&(7.26)\cr
\Xi(u;p,\nu)&=\int_{\B{C}^\times}
j_{0,0}\left(\sqrt{-v}\right)
j_{p,\nu}\left(\sqrt{v/u}\right){d^\times\! v\over|v|}.&(7.27)
}
$$
Then we have
$$
\e{M}(\zeta^2_k,g)=\left\{\e{M}^{(r)}+\e{M}^{(c)}
+\e{M}^{(e)}\right\}(\zeta^2_k,g),\eqno(7.28)
$$ 
with $V$ running over a maximal orthogonal system of
Hecke-invariant cuspidal $\varGamma$-automorphic representations of
$G$. Here $\e{M}^{(r)}(\zeta_k^2,g)$ is analogous to $\e{M}^{(r)}(\zeta^2,g)$, and
$$
\eqalignno{
\e{M}^{(c)}(\zeta^2_k,g)&=\sum_{V}
|\varrho_V(1)|^2H_V\left(\txt{1\over2}\right)^3
\Theta(p_V,\nu_V;g),&(7.29)\cr
\e{M}^{(e)}(\zeta_k,g)&=\sum_{p=-\infty}^\infty\int_{(0)}
{\left|\zeta_k\left({1\over2}(1+\nu),p\right)\right|^6
\over\left|\zeta_k(1+\nu,2p)\right|^2}\Theta(4p,\nu;g){d\nu\over2\pi
i},&(7.30) 
}
$$
where $\zeta_k(s,p)$ is defined by
$$
\zeta_k(s,p)={1\over4}\sum_{n\ne0}(n/|n|)^{4p}|n|^{-2s},\quad \Re
s>1,\eqno(7.31)
$$
which continues meromorphically to $\B{C}$. 
\medskip
The spectral decomposition $(7.28)$ was proved in [8]. 
The argument was a faithful
extension of the older proof of $(6.3)$; that is, it depended on the sum formula of
Kloosterman sums under the situation $(0.3)$ that was established also in [8]. Thus,
in much the same mechanism as Lemmas 1 and 2 did with Kuznetsov's sum
formula,  the last two lemmas should allow us to dispense
with the sum formula of Kloosterman sums, in deriving $(7.28)$.
This claim is still to be checked fully; but it is certain such a new
proof is available and conceptually simpler than that in [8]. 
\smallskip
Comparing $(7.26)$--$(7.30)$ with $(6.1)$--$(6.5)$, the outward similarity
is striking. However, our good luck ends there. That is, the asymptotic
nature of $(6.3)$ is much superior than that of its counterpart $(7.28)$. In fact,
$(7.28)$ does not seem suitable to be utilized as a means to derive 
quantitative assertions on the fourth power moment of
$\zeta_k$. Concerning this difficulty, Sarnak kindly suggested us to try to take
further averaging:
$$
\sum_{q=-\infty}^\infty\e{M}(\zeta_k^2(\cdot,q),g)h(q),\eqno(7.32)
$$
with a smooth weight $h$. The argument of [8] should extend to this sum.
The same can be said about an obvious analogue of $(2.7)$. The latter is
expected to yield an extension of the main result of [19] to 
$L$-functions associated with the group $\r{PSL}_2(\B{Z}[i])$. 
In fact this has been in our current  investigation.
However, we have not worked out all the details yet. 
\bigskip
\noindent
{\bf 8.} Now, we are about to prove Lemmas 3 and 4. This section is a reworking
of  our joint work [7] with Bruggeman; a few corrections are made. It should be
stressed that the surjectivity assertion in Lemma 3 is a new addition, and that the
definition $(7.23)$ is somewhat different from that employed in [32].
\smallskip
We begin with the unitarity  of $\e{K}$:
Naturally it is sufficient to show the orthogonality relation:
$$
{2\over\pi}\int_{\B{C}^\times}\e{K}\phi_{l,q}(u)
\overline{\e{K}\phi_{l',q'}(u)}d^\times\! u
=\delta_{l,l'}\delta_{q,q'},\eqno(8.1)
$$ 
with Kronecker deltas. By definition,
$$
\e{K}\phi_{l,q}(u)=(u/|u|)^{-q}\e{A}_1\phi_{l,q}(\r{a}[|u|])
/\Vert\Phi_{p,q}^l\Vert_K.\eqno(8.2)
$$ 
Note
$$
\eqalignno{
\e{A}_1\phi_{l,q}(\r{g})
&=\exp(2\pi i\Re(z))\sum_{|m|\le l}v_m^l(r)\Phi_{m,q}^l(\r{k})\cr 
&=\exp(2\pi i\Re(z))\sum_{|m|\le
l}e^{-im\varphi-iq\psi}v_m^l(r)\Phi_{m,q}^l(\r{v}[i\theta]), &(8.3)
}
$$ 
where $(7.11)$ has been used, and
$$
\eqalignno{
v^l_m(r)&=\e{A}_1\phi_{l,m}(\r{a}[r])\cr
&=r^{1-\nu}\int_{\B{C}}{\exp(-2\pi ir\Re
v)\over(1+|v|^2)^{1+\nu}}\Phi_{p,m}^l\left(\r{k}\left[
{\overline{v}\over\sqrt{1+|v|^2}},{-1\over\sqrt{1+|v|^2}}\right]
\right)d^+\!v\cr
&=2\pi i^{-p-m}r^{1-\nu}\int_0^\infty{J_{p+m}(2\pi rv)\over(1+v^2)^{1+\nu}}
\Phi_{p,m}^l\left(\r{k}\left[
{v\over\sqrt{1+v^2}},{-1\over\sqrt{1+v^2}}\right]
\right)vdv,&(8.4)
}
$$
where $d^+\!v=(d\Re v)(d\Im v)$.
The left side of $(8.1)$ is equal to
$$
4{\delta_{q,q'}\over\Vert\Phi_{p,q}^l\Vert_K^2}
\int_0^\infty v_q^l(r)\overline{v_q^{l'}(r)} {dr\over r}.\eqno(8.5)
$$ 
\par  
Then we observe that the functions $v_q^l(r)$ satisfy the
differential equations
$$
\eqalignno{ D_q^+ v_q^l(r)&=-4\pi i(l-q)r^{-1}v_{q+1}^l(r),
\cr
D_q^- v_q^l(r)&=4\pi i(l+q)r^{-1}v^l_{q-1}(r), &(8.6)
}
$$ 
where $D_q^-=\overline{D_{-q}^+}$ and
$$ 
D_q^+=\left({d\over dr}\right)^2-(2q+1)r^{-1}{d\over dr}
+r^{-2}(q^2+2q-4\pi^2r^2-8\chi_V^+).\eqno(8.7)
$$ 
To show this, we apply $\Omega_+$ to $(8.3)$. We have
$\Omega_+\e{A}_1\phi_{l,q}(\r{g})=\e{A}_1
\Omega_+\phi_{l,q}(\r{g})$. Thus, by $(7.4)$,
$$
\eqalignno{
\chi_V^+\e{A}_1\phi_{l,q}(\r{g})&=\exp(2\pi i\Re
z)\sum_{|m|\le l}\Bigg\{-{1\over2}\pi^2r^2 +{1\over2}\pi
mre^{i\varphi}\cot\theta+{1\over2}\pi re^{i\varphi}\partial_\theta-\pi
q{re^{i\varphi}\over2\sin\theta}\cr
&+{1\over8}r^2\partial_r^2-{1\over4}mr\partial_r+{1\over8}m^2
-{1\over8}r\partial_r+{1\over4}m\Bigg\}e^{-im\varphi-iq\psi}v_m^l(r)
\Phi_{m,q}^l(\r{v}[i\theta]).\qquad\quad
&(8.8) }
$$
On the other hand, invoking the first line of $(7.12)$ and $(7.14)$, we have
$$
\left( m\cot\theta+
\partial_\theta-{q\over\sin\theta}\right)\Phi_{m,q}^l(\r{v}[i\theta])=
i(l-m+1)\Phi_{m-1,q}^l(\r{v}[i\theta]).\eqno(8.9)
$$
In the last two identities we set $\r{g}=\r{a}[r]$, and note that $\Phi_{m,q}^l(1)
=\delta_{m,q}$. Then we get the first identity of $(8.6)$. In much the same way we 
get the second as well.
\smallskip
Returning to the integral in $(8.5)$, we see that it is equal to
$$
\eqalignno{ 
&-{1\over4\pi i(l-q+1)}\int_0^\infty
D_{q-1}^+v_q^l(r)\cdot
\overline{v_q^{l'}(r)}dr\cr =&-{1\over4\pi i(l-q+1)}
\int_0^\infty v_{q-1}(r)\cdot\overline{D_{q}^-v_q^{l'}(r)}dr\cr
=&{l+q\over l-q+1}\int_0^\infty v_{q-1}^l(r)
\overline{v_{q-1}^{l'}(r)}{dr\over r}. &(8.10)
}
$$ 
This procedure is valid only if $v_q^l(r)$ tends to $0$
sufficiently fast as $r$ tends to either $0$ or $\infty$,
which is in fact implied by the second line of $(8.4)$. Hence
$$
\eqalignno{
\int_0^\infty v_q^l(r)\overline{v_q^{l'}(r)}{dr\over r}
&={l-q\over l+q+1}\int_0^\infty v_{q+1}^l(r)\overline
{v_{q+1}^{l'}(r)}{dr\over r}\cr
&=\delta_{l,l'}\left({2l\atop l-q}\right)^{-1}
\int_0^\infty |v_l^l(r)|^2{dr\over r}.&(8.11)
}
$$ 
On the other hand we have, by the third line of $(8.4)$ and by a formula of N.J. Sonine 
([40, eq.\ (2) on p.\ 434]),
$$
\eqalignno{
v^l_l(r)&=2(-1)^{l-p}i^{-l-p}\pi r^{1-\nu}{2l\choose l-p}\int_0^\infty 
{J_{l+p}(2\pi rv)\over(1+v^2)^{l+1+\nu}}v^{l+p+1}dv\cr
&=2(-1)^{l-p}i^{-l-p}\pi r^{1-\nu}{2l\choose l-p}
{(\pi r)^{l+\nu}\over\Gamma(l+\nu+1)}K_{p-\nu}(2\pi r),  &(8.12)
}
$$
which gives, via either [14, eq.\ 4 of Section 6.576] or [30, $(2.6.11)$],
$$
\int_0^\infty |v_l^l(r)|^2{dr\over r}={1\over4(l+{1\over2})}
\left({2l\atop l-p}\right),\eqno(8.13)
$$ 
and via $(7.10)$, $(8.5)$ and $(8.11)$ we end the proof of $(8.1)$.
\medskip
We turn to the surjectivity assertion. Thus, let $F(u)$, $u\in\B{C}^\times$,
be smooth and compactly supported, and such that 
$$
\int_{\B{C}^\times}F(u)\e{K}\phi_{l,q}(u)d^\times\! u=0,\quad 
\hbox{for all $(q,l)$ with $|p|\le l$, $|q|\le l$}.\eqno(8.14)
$$
We are to show $F\equiv0$. In view of $(8.2)$ we may assume that $F$ is radial, i.e.,
the Fourier expansion of $F$ in $u/|u|$ has only one term, say, the $q$-th. 
With an obvious change in $F$, we consider instead 
$$
\int_0^\infty F(r)v_q^l(r)dr=0,\quad 
\hbox{for all $(q,l)$ with $|p|\le l$, $|q|\le l$}.\eqno(8.15)
$$
We then invoke that in [8, Lemma 5.1] more than $(8.12)$ is proved; thus
there are non-zero $\eta_{p,q}^l(j;\nu)$ such that
$$
v_q^l(r)=\sum_{j=0}^{l-\max\{|p|,|q|\}}\eta_{p,q}^l(j;\nu)
r^{\max\{|p|,|q|\}+1+j}K_{\nu+\max\{|p|,|q|\}-|p+q|+j}(2\pi r).\eqno(8.16)
$$
Namely, we are given, for all integers $l\ge \max\{|p|,|q|\}$,
$$
\int_0^\infty F(r)r^{l+1}K_{\nu+l-|p+q|}(2\pi r)dr=0.\eqno(8.17)
$$
We replace the Bessel factor by a well-known integral representation
(see e.g., [30, (1.1.17)]), and find, after some rearrangement, that $(8.17)$ is the
same as
$$
\int_0^\infty\exp(-\pi \xi)\xi^{\nu+l-|p+q|-1}
\left\{\int_0^\infty F(r)r^{1-\nu+|p+q|}\exp\left(-\pi r^2/\xi\right)dr\right\}
d\xi=0.\eqno(8.18)
$$
Because of the completeness of polynomials over $[0,\infty)$, the member inside
the braces should vanish for any $\xi>0$; that is, for $\Re\xi>0$ by analytic
continuation. Hence, the choice $\xi=1/(1+it)$, $t\in\B{R}$, yields that the Fourier
transform of a multiple of $F(\sqrt{r})$ by a non-zero factor vanishes
constantly. This ends the proof of Lemma 3.
\medskip
 We now move to the proof of Lemma 4. It should be noted that with
ordinary bounds for Bessel functions one may verify absolute
convergence and analytic continuation needed to carry out the reasoning below.
\smallskip 
Thus, we consider the integral
$$
\Gamma_{l,q}(s)=\int_0^\infty v_q^l(r)r^{2(s-1)}dr.
\eqno(8.19)
$$ 
The third line of $(8.4)$ gives
$$
\Gamma_{l,q}(s)=\pi^{1+\nu-2s}(-1)^{\min(0,p+q)}i^{-p-q}
{\Gamma(s+{1\over2}(|p+q|-\nu))\over\Gamma(1-s
+{1\over2}(|p+q|+\nu))}\r{L}_{l,q}(s),\eqno(8.20)
$$ 
with
$$ 
\r{L}_{l,q}(s)=\int_0^\infty{v^{1+\nu-2s}\over
(1+v^2)^{1+\nu}}\Phi^l_{p,q}\left(\r{k}\left[
{v\over\sqrt{1+v^2}},{-1\over\sqrt{1+v^2}}\right]
\right)dv,\eqno(8.21)
$$
since for $m\in\B{Z}$ and $-{1\over2}|m|<\Re s<{1\over4}$
$$
\int_0^\infty J_m(r)(r/2)^{2s-1}dr=(-1)^m
{\Gamma(s+{1\over2}|m|)\over\Gamma(1-s+{1\over2}|m|)}.\eqno(8.22)
$$
We have the functional equation 
$$ 
\r{L}_{l,q}(s)=(-1)^{l-p}\r{L}_{l,-q}(1-s), \eqno(8.23)
$$ 
which is a result of the change of variable
$v\to v^{-1}$ in $(8.21)$. The necessary absolute convergence, and the
meromorphic continuation to $\B{C}$ of $\r{L}_{l,q}(s)$ can be
confirmed readily. Hence we have the local functional equation
$$
\eqalignno{
(-1)^{l-q}&\Gamma_{l,-q}(s)=\pi^{2-4s}(-1)^{\max(|p|,|q|)}
\Gamma_{l,q}(1-s) \cr &\times
{\Gamma(s+{1\over2}(|p+q|+\nu))
\Gamma(s+{1\over2}(|p-q|-\nu))
\over\Gamma(1-s+{1\over2}(|p+q|+\nu))\Gamma(1-s+{1\over2}
(|p-q|-\nu))}&(8.24)
}
$$ 
(see [17, Theorem 6.4]).
\smallskip
Then we observe, by convolving $(8.22)$, that 
$$
\eqalignno{ 
&\int_0^\infty \lambda^{2\nu-1}J_{|p+q|}(r\lambda)
J_{|p-q|}(r/\lambda)d\lambda\cr
&\longleftrightarrow 2^{s-3}
{\Gamma(\txt{1\over4}s+\txt{1\over2}(|p+q|+\nu))
\Gamma(\txt{1\over4}s +\txt{1\over2}(|p-q|-\nu))\over
\Gamma(1-\txt{1\over4}s+\txt{1\over2}(|p+q|+\nu))
\Gamma(1-\txt{1\over4}s +\txt{1\over2}(|p-q|-\nu))}&(8.25)
}
$$ 
is a Mellin pair, provided $2|\Re\nu|<\Re s<1-2|\Re\nu|$. Thus,
denoting the left side by $K_{\nu,p}(r,q)$,
we get, by $(8.19)$, $(8.24)$ and the Mellin--Parseval formula,
$$
(-1)^{l-q}\lambda^{-2}v_{-q}^l(\lambda^2)=8\pi^2(-1)^{\max(|p|,|q|)}
\int_0^\infty
K_{\nu,p}(2\pi \lambda r,q)v_q^l(r^2)rdr.\eqno(8.26)
$$ 
In this we set $\lambda=|u|$ with $u\in\B{C}^\times$, 
and multiply both sides by the factor
$(u/|u|)^{2q}/\Vert\Phi_{p,q}^l\Vert_K$. On noting that 
$\Phi_{p,q}(\r{k}[-\beta,\alpha])=
(-1)^{l-q}\Phi_{p,-q}(\r{k}[\alpha,\beta])$ or $R_{\r{w}}\phi_{l,q}
=(-1)^{l-q}\phi_{l,-q}$, we
see by the definition $(8.2)$ that $(8.26)$ is identical to
$$
\eqalignno{
|u|^{-2}&\e{K}R_\r{w}\phi_{l,q}(u^2) \cr
=&4\pi(-1)^{\max(|p|,|q|)}\int_{\B{C}^\times}
K_{\nu,p}(2\pi |uv|,q)(uv/|uv|)^{2q}\e{K}\phi_{l,q}(v^2)|v|^2
d^\times\!v.\cr
=&4\pi\int_{\B{C}^\times}\left\{\sum_{m=-\infty}^\infty
(-1)^{\max(|p|,|m|)}K_{\nu,p}(2\pi|uv|,m)
\left({uv\over|uv|}\right)^{2m}\right\}
\e{K}\phi_{l,q}(v^2)|v|^2d^\times\! v.\quad&(8.27)
}
$$ 
We then invoke that Graf's addition theorem ([40, eq.\ (1) on p.\
359]) gives, for any $Z,z>0$,
$$
\eqalignno{
\sum_{m=-\infty}^\infty& (-1)^{\max(|p|,|m|)}J_{|m+p|}(Z)
J_{|m-p|}(z)e^{2mi\theta}\cr
&=(-1)^p
J_{2p}\left(|Ze^{i\theta}+ze^{-i\theta}|\right)\left({Ze^{i\theta}+ze^{-i\theta}\over
|Ze^{i\phi}+ze^{-i\theta}|}\right)^{2p}.&(8.28)
}
$$
We apply this to the member inside the braces of $(8.27)$, and find that the proof
of $(7.25)$ with $\phi=\phi_{l,q}$ has been reduced to that of
$$
\eqalignno{
&j_{p,\nu}(u)/(4\pi|u|^2)\cr
&=(-1)^p\int_0^\infty\lambda^{2\nu-1}
J_{2p}\left(2\pi|u||\lambda
e^{i\vartheta}+(\lambda e^{i\vartheta})^{-1}|\right)\left({\lambda
e^{i\vartheta}+(\lambda e^{i\vartheta})^{-1}\over
|\lambda e^{i\vartheta}+(\lambda e^{i\vartheta})^{-1}|}\right)^{2p}d\lambda,
&(8.29)
}
$$
with $u=|u|e^{i\vartheta}$, which is, however, the same as [8, Theorem 12.1].  We end
the proof; the extension to any smooth $\phi$ is immediate.
\medskip
It does not seem that the identity $(8.29)$ had been tabulated before 
[8], a fact somewhat bizarre against its classical appearance. This integral
representation of the Bessel function of representations of
$\r{PSL}_2(\B{C})$ is quite important, for it allows us to deal with test functions
which not necessarily decay exponentially. This merit of $(8.29)$ is indeed exploited
fully in the proof of the spectral decomposition $(7.28)$. The proof in [8] of $(8.29)$
is  conceptually involved, depending for instance on the Goodman--Wallach operator 
([13]). The procedure above indicates the existence of a simpler approach, and
in fact an alternative proof has been obtained in [7][32]. Any extension of $(8.29)$  
is highly desirable.
\smallskip
Lemmas 3 and 4 allow us to carry over the method of [11] to the complex
situation, so that the proof of the Kloosterman--Spectral sum formula established
in [8, Theorem 13.1]  can now be proved in a
simpler manner, although we have not worked out the details yet. 
Finally, we should mention that our argument seems to extend to Lie grpoups of
real rank one; thus, the assertions due to R. Miatello and N.R. Wallach [22] are
hoped to be included in our future discussion.
\vskip 1cm
\noindent
{\bf References}
\medskip
\noindent
\item{[1]}  F.V. Atkinson. The mean value of the Riemann
zeta-function. Acta Math.,\ {\bf 81} (1949), 353--376.
\item{[2]}  E.M. Baruch and  Z. Mao. Bessel identities
in Waldspurger correspondence, the archi\-medean theory.  To appear
in Israel J. Math.
\item{[3]} R.W. Bruggeman.  Fourier Coefficients of
Automorphic Forms. Lecture Notes in Math., {\bf 865},
Springer-Verlag, Berlin etc.\ 1981.
\item{[4]}  R.W. Bruggeman. Automorphic forms.  Banach Center
Publ., {\bf 17} (1985), 31--74.
\item{[5]} R.W. Bruggeman. Sum formula for $\r{SL}_2(\B{C})$. Unpublished, 1995.
\item{[6]}  R.W. Bruggeman and   Y. Motohashi.  Fourth
power moment of Dedekind zeta-functions of real quadratic number
fields with   class number one.  Functiones et Approximatio, {\bf 29}
(2001), 41--79.
\item {[7]}  R.W. Bruggeman  and  Y. Motohashi. A note
on  the mean value of the zeta and $L$-functions.\ XIII. Proc.\ Japan
Acad.,  {\bf 78A} (2002), 87--91.
\item{[8]}  R.W. Bruggeman and  Y. Motohashi.  Sum
formula   for Kloosterman sums and the fourth moment of the Dedekind
zeta-function   over the Gaussian number field.  Functiones et
Approximatio, {\bf 31} (2003),   7--76.
\item{[9]} R.W. Bruggeman and  Y. Motohashi. A new approach to the spectral 
theory of the fourth moment of the Riemann zeta-function. 
J. reine angew.\ Math., {\bf 579} (2005), 75--114.
\item{[10]} D. Bump.  Automorphic Forms on $\r{SL}(3,\B{R})$.
Lecture Notes in Math., vol.\ {\bf 1083},
Springer-Verlag, Berlin etc.\ 1984.
\item{[11]}  J.W. Cogdell and  I. Piatetski-Shapiro. The
Arithmetic and Spectral Analysis of Poin\-car\'e series.
Perspectives in Math., {\bf 13}, Academic Press, San Diego 1990.
\item{[12]}  I.M. Gel'fand, M.I. Graev and  I.I.
Pyatetski-Shapiro.  Representation Theory and Automorphic
Functions. W.B. Saunders Company, Philadelphia 1969.
\item{[13]} R. Goodman and N.R. Wallach. Whittaker vectors and
conical vectors.  J.\ Funct.\ Anal., {\bf39} (1980) 199--279.
\item{[14]}  I.S. Gradshteyn  and  I.M. Ryzhik. Tables
of Integrals, Series, and Products.  Academic Press, San Diego 1979.
\item{[15]}  A. Ivi\'c. Mean Values of the Riemann  
Zeta-Function. Tata IFR Lect.\ Math.\ Phys.\ {\bf 82},
Springer-Verlag, Berlin etc.\ 1991.
\item{[16]} A. Ivi\'c.  The Riemann Zeta-Function.\ Theory and
Applications. Dover Publ., Inc., Mineola, New York 2003.
\item{[17]} H. Jacquet and R.P. Langlands. Automorphic
Forms on $\r{GL}(2)$. Lecture Notes in Math., {\bf 114},
Springer-Verlag, Berlin etc.\ 1970.
\item {[18]}  M. Jutila.  Mean values of Dirichlet series via  
Laplace transforms. In: Analytic Number Theory, Proc.\ 39th
Taniguchi Intern.\ Symp.\ Math., Kyoto 1996, ed.\ Y. Motohashi,
Cambridge Univ.\ Press, Cambridge 1997, pp.\ 169--207.
\item{[19]} M. Jutila and Y. Motohashi. Uniform bound for Hecke $L$-functions.
To appear in Acta Math. 
\item{[20]}  M. Jutila and Y. Motohashi. Uniform bounds for Rankin--Selberg
$L$-functions.  To appear in Proc.\ Workshop on  Multiple Dirichlet Series, Proc.\
Symp.\  Pure Math., AMS.
\item{[21]}  A.A. Kirillov.  On $\infty$-dimensional unitary
representations of the group of second-order matrices with elements
from a locally compact field. Soviet Math.\ Dokl., {\bf 4} (1963),
748--752.
\item{[22]} R. Miatello and N.R. Wallach. Kuznetsov formulas for real rank one
groups. {\it J.\ Funct.\ Anal.}, {\bf93} (1990), 171--206.
\item{[23]}  Y. Motohashi. On $\r{SL}(3,\B{Z})$-Ramanujan sums. Unpublished, 1990.
\item{[24]}  Y. Motohashi.  The fourth power mean of the
Riemann zeta-function. In: Proc.\ Conf.\ Analytic Number Theory,
Amalfi 1989, eds.\ E. Bombieri et al., Univ.\ di Salerno, Salerno
1992, pp.\ 325--344.
\item{[25]}  Y. Motohashi.  An explicit formula for the
fourth   power mean of the Riemann zeta-function. Acta Math., {\bf
170} (1993),   181--220.
\item{[26]}   Y. Motohashi. The mean square of Hecke
$L$-series attached to holomorphic cusp forms. RIMS Kyoto Univ.
Kokyuroku, {\bf 886} (1994), 214--227.
\item{[27]} Y. Motohashi. The binary additive divisor problem.
Ann.\ Sci.\  l'Ecole Norm.\ Sup.\ $4^e$ s\'{e}rie, {\bf 27} (1994),
529--572.
\item {[28]} Y. Motohashi. A relation between the Riemann zeta-function and the
hyperbolic Laplacian. Ann.\ Scuola Norm.\ Sup.\ di Pisa, Sci.\ Fis.\ Mat.,
Ser., IV, 22 (1995), 299--313.
\item {[29]} Y. Motohashi. The mean square of Dedekind zeta-functions of quadratic
number fields. In: Sieve Methods, Exponential Sums, and their Applications
in Number Theory: C. Hooley Festschrift (eds.\ G.R.H. Greaves et al.\ ),
Cambridge Univ.\ Press, Cambridge 1997,  pp.\ 309--324.
\item{[30]}  Y. Motohashi.  Spectral Theory of the Riemann
Zeta-Function. Cambridge Tracts in Math.\ {\bf 127}, Cambridge
Univ.\ Press, Cambridge 1997.
\item{[31]} Y. Motohashi. Addition theorem for Whittaker functions 
and geometric sum formula. Unpublished, 2001; Part II, 2002.
\item{[32]}  Y. Motohashi: Projections of Poincar\'e series into
irreducible subspaces.  Unpublished, 2001.
\item{[33]}  Y. Motohashi. A note on the mean value of the
zeta   and $L$-functions.\ XII. Proc.\ Japan Acad., {\bf 78A} (2002), 36--41.
\item{[34]} Y. Motohashi. A functional equation for the spectral
fourth moment of the modular Hecke $L$-functions. In:
Proc. MPIM-Bonn Special Activity on Analytic Number Theory, 
Bonn 2002, Bonner Math. Schrift., {\bf 130} (2003), 19 pages.
\item{[35]} Y. Motohashi. A vista of mean zeta values. RIMS Kyoto Univ.\ Kokyuroku,
{\bf1319} (2003), 113--124; Part II. ibid., {\bf1384} (2004),129--132.
\item{[36]} Y. Motohashi. A note on the mean value of the zeta and
$L$-functions. XIV. Proc.\ Japan Acad., {\bf 80A} (2004), 28--33. 
\item{[37]} P. Sarnak. Fourth moments of Grossencharakteren 
zeta-functions.\ Comm.\ Pure Appl. Math., {\bf 38} (1985), 167--178.
\item{[38]}  N.Ja.\ Vilenkin.  Special Functions and the
Theory of Group Representations. Amer.\ Math.\ Soc., Providence 1968.
\item{[39]}  N.Ja.\ Vilenkin  and   A.U. Klimyk. 
Representations of Lie Groups and Special Functions.\ Vol.\ 1. Kluwer
Acad.\ Publ., Dordrecht etc.\ 1991.
\item{[40]}  G.N. Watson. A Treatise on the Theory of Bessel
Functions. Cambridge Univ.\ Press, Cambridge 1996.
\item{[41]}  E.T. Whittaker and  G.N. Watson.  A
Course of Modern Analysis. Cambridge Univ.\ Press, London 1969.
\vskip 1cm
\font\small=cmr8 
\noindent {\small Yoichi Motohashi
\par\noindent Department of Mathematics, Nihon University,
\par\noindent Surugadai, Tokyo 101-8308, Japan
\par\noindent Email: ymoto@math.cst.nihon-u.ac.jp}

\bye